\documentclass[12pt,a4paper]{article}
\usepackage{amsfonts, amsmath, amssymb,latexsym}
\usepackage{theorem}
\usepackage[english]{babel}
\usepackage{epsfig}
\usepackage{enumerate}
\usepackage{array}

\usepackage[latin1]{inputenc}

\newtheorem{thm}{Theorem}[section]
\numberwithin{equation}{section}
\newtheorem{teorema}[thm]{Theorem}
\newtheorem{prop}[thm]{Proposition}
\newtheorem{coro}[thm]{Corollary}

\theorembodyfont{\rm}

\newtheorem{Remarque}[thm]{Remark}

\newcommand{\SL}{\mathrm{SL}(2,\mathbb{Z})}

\newcommand{\PGL}{\mathrm{PGL}(2,\mathbb{Z})}
\newcommand{\PSL}{\mathrm{PSL}(2,\mathbb{Z})}

\newcommand{\A}{\mathcal{A}}
\newcommand{\B}{\mathcal{B}}

\newcommand{\Red}{\mathrm{Red}}
\newcommand{\Sim}{\mathrm{Sim}}

\newcommand{\discr}{\mathrm{discr}}

\def\={\;=\;}
\def\.={\;\dot{=}\;}

\def\O{\mathcal{O}}

\usepackage{fancyhdr}

\begin{document}
\title{From continued fractions and quadratic functions to modular forms}
\author{Paloma Bengoechea}

\maketitle

\begin{abstract} In this paper we study certain real functions defined in a very simple way by Zagier \cite{Z99} as sums of infinite powers of quadratic polynomials with integer coefficients. These functions give the even parts of the period polynomials of the modular forms which are the coefficients in Fourier expansion of the kernel function for Shimura-Shintani correspondence.
We prove two conjectures of Zagier showing that the sums converge exponentially. We also prove unexpected results on the representation of these functions as sums over simple or reduced quadratic forms and the positive or negative continued fraction of the variable. These arise from more general results on polynomials of even degree. Especially we give the even part of the Eichler integral on $x\in\mathbb{R}$ of any cusp form for $\PSL$ in terms of the even part of its period polynomial and the continued fraction of $x$.
 
\end{abstract}

\section{Introduction}

In \cite{Z99}, D. Zagier studied certain functions defined as sums of powers of quadratic polynomials with integer coefficients and discovered that these functions have many surprising properties. For example, for $x\in\mathbb{R}$, take the sum over all quadratic functions $Q(X)=aX^2+bX+c$ with integer coefficients and fixed discriminant that are negative at infinity and positive at $x$. Then one finds (for convenience we write ``$Q(\infty)<0$'' to mean that $a_Q=a$ is negative):

\begin{teorema}\label{Zagier} Let $D$ be a positive non-square discriminant. Then the sum
\begin{equation}\label{A}
A_D(x):=\sum_{\substack{\discr(Q)=D\\Q(\infty)<0<Q(x)}} Q(x)
\end{equation}
converges for all $x\in\mathbb{R}$ and has a constant value $\alpha_D$ independent of $x$.
\end{teorema}

For example, one finds $A_5(x)=2$, and more generally $\alpha_D=-5L(-1,\chi_D)$, where $L(s,\chi_D)$ is the Dirichlet $L$-series of the character $\chi_D(n)=(D/n)$ (Kronecker symbol). We denote
$$
\mathcal{Q}=\left\{aX^2+bX+c\mid\, a,b,c\in\mathbb{Z},\, b^2-4ac>0,\, b^2-4ac\mbox{ not a square}\right\},
$$
and for any positive non-square discriminant $D$,
\begin{align*}
&\mathcal{Q}_D=\left\{aX^2+bX+c\in\mathcal{Q}\mid\, b^2-4ac=D\right\},\\
&\mathcal{Q}_D\langle x\rangle=\left\{Q\in\mathcal{Q}_D\mid\, Q(\infty)<0<Q(x)\right\}.
\end{align*} 
The sum in \eqref{A} is the special case $k=2$ of the function
\begin{equation}\label{A k,D}
A_{k,D}(x):=\sum_{Q\in\mathcal{Q}_D\langle x\rangle} Q(x)^{k-1}\qquad (k\geq 2).
\end{equation}
Theorem \ref{Zagier} is still true for $k=4$ with $\alpha_D$ now replaced by $L(-3,\chi_D)$, but this fails for even $k\geq 6$ because of the existence of cusp forms of weight $2k$ on the full modular group. More explicitly, for even $k\geq 6$, the function $A_{k,D}(x)$ is a linear combination of a constant function and the functions $\sum_{n\geq 1}n^{1-4k}a_f(n)\cos(2\pi inx)$, where $f$ runs over the normalized Hecke eigenforms in $S_{2k}(\PSL)$ and $a_f(n)$ denotes the $n$-th Fourier coefficient of $f$.
The function $A_{k,D}(x)$ arose from studying the modular form of weight $2k$ which is the $D$-th coefficient in Fourier expansion of the kernel function for the Shimura and Shintani lifts between half-integral and integral weight cusp forms (\cite{Shim73},\cite{Shin75}). The function $A_{k,D}(x)$ is (modulo a constant multiple of $X^{2k-2}-1$) the even part of its Eichler integral on $\mathbb{R}$ and so gives the even part of its period polynomial (\cite{KZ80}, \cite{KZ84} and \cite{Z99}).

The convergence of $A_{k,D}(x)$ is not immediate at all. As Zagier observed, one has $Q(x)=O(1/a_Q)$ for all the $Q$ occuring in \eqref{A k,D} and one sees easily that there are only O(1) quadratic functions $Q$ for each $a$-value. Therefore the series \eqref{A k,D} converges at most like $\sum_{a>0}a^{1-k}$ if $k\geq 4$. 
In the case $k=2$, however, this argument fails and Zagier could only deduce the convergence of \eqref{A} from the fact that the function is finite and constant ($=-5L(-1,\chi_D)$) when $x$ is rational. In fact, if for some value of $x$ the sum diverged, then, since all summands are positive, there would be finitely many quadratic functions $Q$ whose sum at $x$ already exceeded $-5L(-1,\chi_D)$, and the sum of their values at a sufficiently nearby rational number would also exceed $-5L(-1,\chi_D)$.

Following this proof, the sum might converge extremely slowly.
However, experiments carried out for the value $x=1/\pi$ suggested that the series \eqref{A} (and hence also the series \eqref{A k,D}) converge extremely rapidly. More precisely, for $x=1/\pi$ and $D=5$, Zagier found experimentally that the elements of $\mathcal{Q}_D\langle x\rangle$ belong to the union of two lists, each having values that tend to 0 exponentially quickly. 
Each quadratic function $Q$ in each list is obtained from the preceding one by applying a fairly simple element of $\PGL$ giving the positive continued fraction of $x$. The first functions in the lists correspond to the opposite of the simple forms coming from the reduction theory of binary quadratic forms with fixed discriminant 5: $[1,-1,-1]$ and $[1,1,-1]$.
The following tables give the first five functions $Q$, and the corresponding values of $Q(1/\pi)$ for each list:
\begin{center}
\begin{tabular}{c|c}
 $Q$ &$Q(1/\pi)$ \\
\hline
$\left[-1,1,1\right]$ & 1.216989 \\ 
$\left[-11,7,-1\right]$ & 0.113636 \\ 
$\left[-541,345,-55\right]$  &0.00215\\ 
$\left[-117731,74951,-11929\right]$ &0.000008  \\ 
$\left[-133351,84893,-13511\right]$  &0.000008  \\ \hline
Sum: &1.332791
\end{tabular}
\end{center}

\begin{center}
\begin{tabular}{c|c}
$Q$ & $Q(1/\pi)$ \\
\hline
$\left[-1,-1,1\right]$ &0.580369\\ 
$\left[-5,5,-1\right]$ &0.084943\\ 
$\left[-409,259,-41\right]$  &0.001896\\ 
$\left[-5959340757998441,3793834156817819,-603807459328429\right]$ &6.856501 E-17\\ 
$\left[-7755390254828071,493723477865040,-785785320227431\right]$ &1.568047 E-18\\ \hline
Sum: &0.667208
\end{tabular}
\end{center}

Zagier conjectured that all functions in one list and some in the second one occur in $A_5(x)$, but he found no criterion to decide which ones. 
He suggested there is a similar situation in the general case. 

We will prove Zagier's conjectures, giving a criterion for the second list, and establish a similar result for the general case, obtaining the exponential convergence for $A_{k,D}(x)$ and a direct proof for the convergence in the case $k=2$. 
We also prove other descriptions of $A_{k,D}(x)$ analogous to the one conjectured by Zagier: first we replace the simple forms with the reduced forms and the usual algorithm of positive continued fraction with a different one (with sign + as well). Later we keep the simple forms but we consider the (usual) negative continued fraction. We give correspondences between all these situations.
This is done in the third section.

In the fourth section we prove an unexpected result: the values at $x$ of the functions in the lists that do not appear in the sum $A_{k,D}(x)$ cancel each other out. Thus the sum over all functions in the lists is equal to $A_{k,D}(x)$. A similar phenomenon holds for the description with the negative continued fraction. This is a consequence of the more general Corollary \ref{prop funcion auxiliar} which gives the even part of the Eichler integral on $x\in\mathbb{R}$ of a cusp form for $\PSL$ in terms of the even part of its period polynomial and the continued fraction of $x$.

\section{Reduction theories}

In this section we give the main connections between reduction theory of binary quadratic forms with positive non square discriminant for $\PSL$ and the continued fractions of a real number. We also recall some simple properties of positive and negative continued fractions and give the non usual algorithm that we use in the next sections.

We denote 
$$
\varepsilon=\begin{pmatrix} 0 &1\\1 &0\end{pmatrix},\qquad
\sigma=\begin{pmatrix} -1 &0\\ 0 &1\end{pmatrix},\qquad
S=\begin{pmatrix} 0 &-1\\ 1 &0\end{pmatrix},\qquad
T=\begin{pmatrix} 1 &1\\ 0 &1\end{pmatrix}.
$$
The matrices $\varepsilon$,$\sigma$,$T$ generate the group $\Gamma=\PGL$ and $S$,$T$ generate the group $\Gamma_1=\PSL$. Clearly 
$$
\varepsilon^2=\sigma^2=S^2=\varepsilon\sigma S=1
$$
so $\left\{1,\varepsilon,\sigma,S\right\}$ form a Klein 4-group. We write $\hat{\Gamma}$ to mean $\Gamma$ or $\Gamma_1$. The group $\hat{\Gamma}$ acts on the set of binary quadratic forms by
\begin{equation}\label{accion}
Q|\gamma(X,Y)=Q(rX+sY,tX+uY)\qquad (\gamma=\begin{pmatrix}r &s\\t &u\end{pmatrix}\in\hat{\Gamma}).
\end{equation}

For a positive non-square fixed discriminant $D$, a \textit{reduction theory} of binary quadratic forms with discriminant $D$ for the action of the group $\Gamma_1$ consists in giving a finite system $\mathcal{R}$ of such forms (called a system of reduced forms) and an algorithm such that:

(i) each form with discriminant $D$ is $\Gamma_1$-equivalent to some element of $\mathcal{R}$ by applying the algorithm a finite number of times;

(ii) the image by the algorithm of an element of $\mathcal{R}$ still belongs to $\mathcal{R}$. In other words, the elements of $\mathcal{R}$ form cycles such that two elements of $\mathcal{R}$ are $\Gamma_1$-equivalent if and only if they belong to the same cycle. In particular, the number of cycles is the number of $\Gamma_1$-equivalence classes for $D$.

Usually (and in all cases we consider) the reduction algorithm for binary quadratic forms $Q=[a,b,c]$ ($:=aX^2+bXY+cY^2$) is obtained by applying a reduction algorithm for real numbers (usually some sort of continued fraction algorithm) to one of the roots of $Q(X,1)$ 
$$
w_Q=\frac{-b-\sqrt{D}}{2a},\qquad w_Q'=\frac{-b+\sqrt{D}}{2a},
$$
where $\sqrt{D}$ denotes the positive square root.


We use two different sets $\mathcal{R}$ of forms $[a,b,c]$:
$$
 a>0,\qquad c>0,\qquad b>a+c
$$
which we call \textit{reduced} (see \cite{Z81} for details), and  
$$
 a>0>c
$$
which we call \textit{simple}.

The bijection $[a,b,c]\mapsto[-a,-b,-c]$ exchanges the simple forms with positive and negative values of $a+b+c$ (the value 0 cannot occur  for non-square $D$).
The bijection
$$
\begin{array}{lll}
\mbox{reduced} &\rightarrow &\mbox{simple with}\, \, a+b+c>0\\
\left[a,b,c\right] &\mapsto &\left[a,b-2a,c-b+a\right]
\end{array}
$$
proves that there are exactly half as many reduced forms as simple forms.
We denote by $\Red$ and $\Sim$ the sets of quadratic polynomials that correspond to the reduced and simple forms respectively
$$
\Red=\left\{Q(X)\in\mathcal{Q}\, |\, Y^2Q(X/Y)\, \mbox{is reduced}\right\},
$$
$$
\Sim=\left\{Q(X)\in\mathcal{Q}\, |\, Y^2Q(X/Y)\, \mbox{is simple}\right\}.
$$
We can translate the inequalities for the coefficients of simple or reduced forms into inequalities for the corresponding quadratic irrationalities:
$$
\begin{array}{rl}
Q(X,Y)\,\, \mbox{reduced}\quad &\Longleftrightarrow\quad w_Q<-1<w_Q'<0,\\\\
Q(X,Y)\,\, \mbox{simple}\quad &\Longleftrightarrow\quad w_Q<0<w_Q'.
\end{array}
$$

The positive and negative continued fraction of a real number $x$, denoted by
$$
x=n_0+\dfrac{1}{n_1+\dfrac{1}{n_2+\dfrac{1}{\ddots}}}\qquad (n_i\in\mathbb{Z},\, n_i\geq 1\ \forall i\geq 1),
$$
$$
x=m_0-\dfrac{1}{m_1-\dfrac{1}{m_2-\dfrac{1}{\ddots}}}\qquad(m_i\in\mathbb{Z},\, m_i\geq 2\ \forall i\geq 1),
$$
are also produced by reduction algorithms:
\begin{equation}\label{algo Red 1}
x_0=x,\qquad n_i=\left\lfloor x_i\right\rfloor,\qquad x_{i+1}=\dfrac{1}{x_i-n_i}=\varepsilon T^{-n_i}(x_i)\qquad (i\geq 0),
\end{equation}
\begin{equation}\label{algo Red 2}
x_0=x,\qquad m_i=\left\lceil  x_i\right\rceil+1,\qquad x_{i+1}=\dfrac{1}{m_i-x_i}=ST^{-m_i}(x_i)\qquad (i\geq 0),
\end{equation}
where $\left\lfloor x\right\rfloor$ and $\left\lceil x\right\rceil$ are respectively the floor and ceil parts.
The positive and negative continued fractions of $w_Q$ are periodic for a quadratic form $Q$. Moreover, the negative continued fraction of $w_Q$ is purely periodic if and only if $Q$ is reduced.
Thus the cycle of reduced forms that are equivalent to a given form $Q$ corresponds to the cycle of real quadratic irrationalities $x_i$ given by the algorithm \eqref{algo Red 2} with $x_0=w_Q$.

In a similar way, a quadratic form $Q$ is simple if and only if $w_Q$ is purely periodic for the algorithm (\cite{ChZ})
$$
x_0=x,\qquad x_{i+1}=\left\{\begin{array}{llll}
x_i+1 &= &T(x_i) &\mbox{if $x_i\leq 0$},\\
\\
\dfrac{x_i}{1-x_i}&=&T^{-1}ST^{-1}(x_i) &\mbox{if $0<x_i<1$},\\
\\
x_i-1&=&T^{-1}(x_i) &\mbox{if $x_i\geq1$}.
\end{array}\right.
$$
This gives an expansion in negative continued fraction that is slower than the usual one.

The cycle of simple forms which are equivalent to a given form $Q$ corresponds to the cycle of real quadratic irrationalities $x_i$ given by the algorithm above with $x_0=w_Q$.
\\


Clearly each $x_i$ in \eqref{algo Red 1} is the image of $x$ by a matrix $\gamma_i=\gamma_{i,x}\in\Gamma$, given explicitly by 
\begin{equation}\label{matrices primeras}
\gamma_0=\gamma_{0,x}:=\mathrm{Id},\qquad \gamma_i=\gamma_{i,x}:=\begin{pmatrix} 0 &1\\1 &-n_{i-1}\end{pmatrix}\cdots\begin{pmatrix} 0 &1\\1 &-n_0\end{pmatrix}\qquad(i\geq 1)
\end{equation}
and recursively by
\begin{equation}\label{matrices}
\gamma_{0}=\mathrm{Id},\qquad \gamma_{i+1}=\varepsilon T^{-n_i}\gamma_{i}\qquad(i\geq 0).
\end{equation}
We denote 
$$\Gamma(x)\, :=\, \left\{\gamma_1,\gamma_2,\gamma_3,\ldots\right\}\, \subset\, \Gamma.
$$ 
There is an explicit description of $\Gamma(x)$ in terms of the convergents of $x$.
The $i$-th convergent of $x$ is denoted by $\dfrac{p_i}{q_i}=[n_0,\ldots,n_i]$.
The integers $p_i$ and $q_i$ satisfy the recurrence 
$$
\begin{array}{ccccc}
&p_{-2}=0 &\quad p_{-1}=1,&\qquad p_i=n_ip_{i-1}+p_{i-2} &\qquad (i\geq 0),\\
&q_{-2}=1,&\quad q_{-1}=0, &\qquad q_i=n_iq_{i-1}+q_{i-2} &\qquad (i\geq 0),
\end{array}
$$
the equation
\begin{equation}
p_{i+1}q_i-p_iq_{i+1}=(-1)^i
\end{equation}
and the inequalities
\begin{enumerate}[(1)]
\item $q_i\geq q_{i-1}\geq0$ for all $i\geq 0$ and $q_i>q_{i-1}>0$ for all $i\geq 2$,

\item $|p_i|\geq|p_{i-1}|$ for all $i\geq2$ and $|p_i|>|p_{i-1}|$ for all $i\geq3$.


\end{enumerate}
The numbers $\delta_i$ ($i\geq -1$) defined by 
\begin{equation}\label{delta}
\delta_{i}=(-1)^{i}(p_{i-1}-q_{i-1}x)
\end{equation}
satisfy the recurrence
$$
\delta_{-1}=x,\qquad \delta_{0}=1,\qquad \delta_{i+1}=-n_i\delta_{i}+\delta_{i-1}\qquad\mbox{with }n_i=\left\lfloor \dfrac{\delta_{i-1}}{\delta_{i}}\right\rfloor
$$
and the inequalities $1=\delta_{0}>\delta_1>\ldots\geq0$. If $x$ is rational, then $x_i=p_i/q_i$ for some $i$ and the recurrence stops with $\delta_{i+1}=0$; if $x$ is irrational, the $\delta_i$ are all positive and converge to 0 with exponential rapidity. 
With these notations, one has
\begin{equation}
\gamma_i^{-1}=\begin{pmatrix}p_{i-1} &p_{i-2}\\q_{i-1} &q_{i-2}\end{pmatrix},\qquad \gamma_{i}\begin{pmatrix}x\\ 1\end{pmatrix}=\begin{pmatrix}\delta_{i-1}\\\delta_i\end{pmatrix}.
\end{equation}

Now we consider the slower version of the algorithm of reduction \eqref{algo Red 1}
\begin{equation}\label{algo 3}
x_0=x,\qquad x_{i+1}=\left\{\begin{array}{llll}
x_i+1 &= &T(x_i) &\mbox{if $x_i\leq0$},\\
\\
\dfrac{1}{x_i}-1 &= &T^{-1}\varepsilon(x_i) &\mbox{if $0<x_i\leq1$},\\
\\
x_i-1 &= &T^{-1}(x_i) &\mbox{if $x_i>1$},
\end{array}\right.
\end{equation}
such that with this algorithm the expansion of $x$ in continued fraction is
$$
x=\underbrace{\pm1\pm\cdots\pm1}_{|n_0|}+\dfrac{1}{\underbrace{1+\cdots+1}_{n_1}+\dfrac{1}{\underbrace{1+\cdots+1}_{n_2}+\dfrac{1}{\ddots}}}
$$
where $n_0,n_1,n_2,\ldots,$ are given in \eqref{algo Red 1} and each $\pm$ equals the sign of $|n_0|$.

Each $x_i$ in algorithm \eqref{algo 3} is the image of $x$ by a matrix $\gamma_i'=\gamma_{i,x}'\in\Gamma$ given recursively by
\begin{equation}\label{matrices 2}
\gamma'_{0}=\mathrm{Id},\qquad \gamma'_{i+1}=\left\{\begin{array}{ll}
T\gamma'_{i} &\mbox{if $x_i\leq0$},\\\\
T^{-1}\varepsilon\gamma'_{i} &\mbox{if $0<x_i\leq1$},\\\\
T^{-1}\gamma'_{i} &\mbox{if $x_i>1$}.
\end{array}\right.
\end{equation}
We denote
$$
\Gamma(x)'\, :=\, \left\{\gamma_1',\gamma_2',\gamma_3',\ldots\right\}\, \subset\, \Gamma.
$$
We note that 
$$
\Gamma(x)'=\left\{T^{-k}\gamma_{i},\ 1\leq k\leq n_{i}\right\}_{i\geq 1}.
$$ 

The following two propositions, whose proofs are given in \cite{B}, describe the sets $\Gamma(x)$ and $\Gamma(x)'$ for $x\in\mathbb{R}-\mathbb{Q}$ as subspaces of elements of $\Gamma$ defined by certain simple linear inequalities:

\begin{prop}\label{prop 1}
For all $x\in\mathbb{R}$ irrational, the set $\Gamma(x)$ equals $W-(W_1\cup W_2)$, where
$$
W=\left\{\gamma\in\Gamma\, \mid -1\leq\gamma(\infty)\leq 0,\, \gamma(x)>1\right\},
$$
$$
W_1=\left\{\gamma\in W\, \mid\, \gamma(\infty)=0,\, \det(\gamma)=1\right\}\=\left\{\begin{pmatrix}0 &-1\\1 &-1-n_0\end{pmatrix}\right\},
$$
$$
W_2=\left\{\gamma\in W\, \mid\, \gamma(\infty)=-1,\, \mathrm{det}(\gamma)=-1\right\}\=\left\{\begin{array}{cl}\left\{
\begin{pmatrix}-1 &1+n_0\\1 &-n_0\end{pmatrix}\right\} &\quad\mbox{if $n_1\geq2$},\\
\emptyset &\quad\mbox{if $n_1=1$}.
\end{array}\right.
$$
\end{prop}

\begin{Remarque}Proposition \ref{prop 1} is also true for $x\in\mathbb{Q}$ if we allow the value $\infty$ for $\gamma(x)$, when $\gamma\in W$.
\end{Remarque}

\begin{prop} 
For all $x\in\mathbb{R}$ irrational, the set $\Gamma(x)'$ equals $W'-W'_1$, where
$$
W'=\left\{\gamma\in\Gamma\, \mid\, \gamma(\infty)\leq -1,\, \gamma(x)>0\right\}
$$
and
$$
W'_1=\left\{\gamma\in W'\, \mid\, \gamma(\infty)=-1,\, \det(\gamma)=1\right\}\=\left\{\begin{pmatrix}1 &-n_0\\-1 &n_0+1\end{pmatrix}\right\}.
$$
\end{prop}

Each $x_{i}$ in \eqref{algo Red 2} is the image of $x$ by a matrix $\tilde\gamma_i=\tilde\gamma_{i,x}\in\Gamma_1$ defined by 
\begin{equation}\label{matrices 3}
\tilde\gamma_{0}=\mathrm{Id},\qquad \tilde\gamma_{i+1}=ST^{-m_i}\tilde\gamma_{i}=\left(\begin{array}{cc}
-\tilde q_{i-1} &\tilde p_{i-1}\\
-\tilde q_i &\tilde p_i
\end{array}\right) \qquad (i\geq 0),
\end{equation}
where $\dfrac{\tilde p_i}{\tilde q_i}=[m_0,\ldots,m_i]$ is the $i$-th convergent of $x$.
The set of matrices from \eqref{matrices 3} will be denoted by
$$
\Gamma_1(x)\, :=\, \left\{\tilde\gamma_1,\tilde\gamma_2,\tilde\gamma_3,\ldots\right\}\, \subset\, \Gamma_1.
$$
 
The integers $\tilde{p}_i$ and $\tilde{q}_i$ satisfy the recurrence 
$$
\begin{array}{ccccc}
\tilde{p}_{-2}=0,&\qquad \tilde{p}_{-1}=1, &\qquad \tilde{p}_i=m_i\tilde{p}_{i-1}-\tilde{p}_{i-2} &\qquad (i\geq 0),\\
\tilde{q}_{-2}=-1,&\qquad \tilde{q}_{-1}=0, &\qquad \tilde{q}_i=m_i\tilde{q}_{i-1}-\tilde{q}_{i-2} &\qquad (i\geq 0),
\end{array}
$$
the equation
\begin{equation}
\tilde{p}_{i}\tilde{q}_{i+1}-\tilde{p}_{i+1}\tilde{q}_{i}=1,
\end{equation}
and the inequalities
\begin{enumerate}[(1)]
\item $\tilde{q}_i\geq \tilde{q}_{i-1}\geq0$ for all $i\geq 0$ and $\tilde{q}_i>\tilde{q}_{i-1}>0$ for all $i\geq 1$,

\item $|\tilde{p}_i|\geq|\tilde{p}_{i-1}|$ for all $i\geq1$ and $|\tilde{p}_i|>|\tilde{p}_{i-1}|$ for all $i\geq2$.


\end{enumerate}

In a similar way to the positive continued fraction, the numbers $\tilde\delta_i$ ($i\geq -1$) defined by 
\begin{equation}\label{delta}
\tilde\delta_{i}=\tilde p_{i-1}-\tilde q_{i-1}x
\end{equation}
satisfy the recurrence
$$
\tilde\delta_{-1}=x,\qquad \tilde\delta_{0}=1,\qquad\tilde\delta_{i+1}=m_i\tilde\delta_{i}-\tilde\delta_{i-1}\qquad\mbox{with }m_i=\left\lceil \frac{\tilde\delta_{i-1}}{\tilde\delta_{i}}\right\rceil
$$
and the inequalities $1=\tilde\delta_{0}>\tilde\delta_1>\ldots\geq0$. If $x$ is rational, then $x_i=\tilde p_i/\tilde q_i$ for some $i$ and the recurrence stops with $\tilde\delta_{i+1}=0$; if $x$ is irrational, the $\tilde\delta_i$ are all positive and converge to 0 with exponential rapidity. 

\section{Combining reduction theories: proofs of the conjectures}

For $d\in 2\mathbb{N}$, the group $\hat{\Gamma}$ acts on homogeneous polynomials of degree $d$ by \eqref{accion} or, equivalently, on the space of polynomials of degree $\leq d$ in one variable by
$$
(P|_{-d}\gamma)(x):=(tx+u)^d\, P\Big(\dfrac{rx+s}{tx+u}\Big)\qquad (\gamma=\begin{pmatrix}r &s\\t &u\end{pmatrix}\in\hat{\Gamma}).
$$
We write
$$
\mathcal{F}_d=\left\{P\in\mathbb{Z}[X]_{\leq d}\mid\, P \mbox{ has exactly 2 real roots, both irrational}\right\}.
$$
Given $P\in\mathcal{F}_d$, we denote by $w_P$ and $w_P'$ its two real roots such that 
$$
\mathrm{sign}(P(\infty))\cdot w_P<\mathrm{sign}(P(\infty))\cdot w_P'.
$$ 
If $\A$ is a $\hat{\Gamma}$-equivalence class in $\mathcal{F}_d$, we define
\begin{align*} 
&\A^\Red=\left\{P\in\A\mid\, P(\infty)>0,\, P(-1)<0<P(0)\right\}=\left\{P\in\A\mid\, w_P<-1<w_P'<0\right\},\\
&\A^\Sim=\left\{P\in\A\mid\, P(0)<0<P(\infty)\right\}=\left\{P\in\A\mid\, w_P<0<w_P'\right\},\\
&\A\langle x\rangle=\left\{P\in\A\mid\, P(\infty)<0<P(x)\right\}.
\end{align*}

In the special case $d=2$, we have $\mathcal{F}_2=\mathcal{Q}$
and 
$$
\A^{\Red}=\A\cap\Red,\qquad \A^{\Sim}=\A\cap\Sim.
$$
In this case, $\A^\Red$ and $\A^\Sim$ are both finite. If $\A$ is a $\hat{\Gamma}$-equivalence class in $\mathcal{Q}$, we define
$$
A_{k,\A}(x)=\sum_{Q\in\A\langle x\rangle} Q(x)^{k-1}\qquad(x\in\mathbb{R},\, k\geq2).
$$
Then the sum $A_{k,D}$ defined in \eqref{A k,D} is given by
\begin{equation}
A_{k,D}(x)\=\sum_{\A\in\mathcal{Q}_D/\hat{\Gamma}} A_{k,\A}(x).
\end{equation}

The two theorems below are stated in the general case $\mathcal{F}_d$ with even $d\geq 2$.
\begin{teorema}\label{teorema 1}
For $d\in 2\mathbb{N}$, $\A$ a $\Gamma$-equivalence class in $\mathcal{F}_d$ and $x\in\mathbb{R}$, the following bijection holds 
$$
\begin{array}{rll}
\left\{\begin{array}{c}
(P,\gamma)\in\A^\Sim\times\Gamma(x):\\ P(\gamma(\infty))<0<P(\left\lfloor \gamma(x)\right\rfloor)\end{array}
\right\} &\stackrel{\cong}{\longrightarrow} &\A\langle x\rangle\\
(P,\gamma) &\mapsto &P|\gamma.
\end{array}
$$
\end{teorema}

\textbf{Proof.} We only have to check that $P(\left\lfloor \gamma(x)\right\rfloor)>0$ implies $P(\gamma(x))>0$ to prove that the map is well defined. This follows from $P(\left\lfloor \gamma(x)\right\rfloor)>0$ and $\left\lfloor \gamma(x)\right\rfloor>0$, which imply $\left\lfloor \gamma(x)\right\rfloor>w'_P$, so $\gamma(x)>w'_P$, and hence $P(\gamma(x))>0$.

We now prove that the map is a bijection. 
Let $P\in\A$ satisfy $P(\infty)<0<P(x)$. 
For $j\gg0$, the convergents $\dfrac{p_j}{q_j}$ of $x$ belong to $(w_P',w_P)$, because $x\in(w_P',w_P)$. Moreover, if two consecutive convergents $\dfrac{p_j}{q_j}$ and $\dfrac{p_{j+1}}{q_{j+1}}$ belong to $(w_P',w_P)$, so do all the later convergents.

If $\left\lfloor x\right\rfloor\not\in(w_P',w_P)$, we define $i$ to be the unique positive integer such that $\dfrac{p_{i-1}}{q_{i-1}}\not\in(w_P',w_P)$ but $\dfrac{p_i}{q_i},\dfrac{p_{i+1}}{q_{i+1}}\in(w_P',w_P)$. If $\left\lfloor x\right\rfloor\in(w_P',w_P)$, we set $i=0$. Since $P(\infty)<0$, in both cases we have
\begin{equation}\label{unicidad i}
P\Big(\dfrac{p_{i-1}}{q_{i-1}}\Big)<0<P\Big(\dfrac{p_i}{q_i}\Big),\qquad P\Big(\dfrac{p_{i+1}}{q_{i+1}}\Big)>0.
\end{equation}
Put 
$$
\gamma=\begin{pmatrix}
q_{i-1} &-p_{i-1}\\
-q_{i} &p_{i}
\end{pmatrix},
$$
and $R=P|\gamma^{-1}$. We have
$$
R(0)=P\Big(\dfrac{p_{i-1}}{q_{i-1}}\Big)<0,\qquad\qquad R(\infty)=P\Big(\dfrac{p_i}{q_i}\Big)>0.
$$ 
The inequality $P\Big(\dfrac{p_{i+1}}{q_{i+1}}\Big)>0$ is equivalent to the condition $R(\left\lfloor \gamma(x)\right\rfloor)>0$ because $P\Big(\dfrac{p_{i+1}}{q_{i+1}}\Big)=R(n_{i+1})$ and $n_{i+1}=\left\lfloor \gamma(x)\right\rfloor$.
The condition $P(\infty)<0$ is equivalent to the condition $R(\gamma(\infty))<0$. Thus $(R,\gamma)$ belongs to the left hand set of the map in Theorem \ref{teorema 1}.

The uniqueness of the preimage $(R,\gamma)$ follows from the equivalence between the condition $R(\left\lfloor \gamma(x)\right\rfloor)>0$ and the inequality $P\Big(\dfrac{p_{i+1}}{q_{i+1}}\Big)>0$, together with the uniqueness of $i$ satisfying \eqref{unicidad i}.
\begin{flushright}
$\square$
\end{flushright}

Note that we did not use the fact that $\Gamma(x)$ comes from the continued fraction of $x$ to prove that the map above is well defined, but rather the description given by Proposition \ref{prop 1}. The argument for the bijectivity is in fact a ``local'' phenomenon: we did not need the whole continued fraction of $x$, but only three consecutive convergents. One could certainly prove Theorem \ref{teorema 1} without using continued fractions and using the description for $\Gamma(x)$ with linear inequalities, but the proof given here seemed to the author to be simple and attractive.

\begin{Remarque} If we replace $\left\lfloor \gamma(x)\right\rfloor$ by $\gamma(x)$ in the above definition, then the map
$$
\begin{array}{rll}
\left\{\begin{array}{c}
(P,\gamma)\in\A^{\Sim}\times\Gamma(x):\\ P(\gamma(\infty))<0<P(\gamma(x))\end{array}
\right\} &\longrightarrow &\A\langle x\rangle\\
(P,\gamma) &\mapsto &P|\gamma
\end{array}
$$
is still surjective but in general not injective.
\end{Remarque}

\begin{teorema}\label{teorema 2}
For $d\in 2\mathbb{N}$, $\A$ a $\Gamma$-equivalence class in $\mathcal{F}_d$ and $x\in\mathbb{R}$, the following bijection holds
$$
\begin{array}{rll}
\left\{(P,\gamma)\in\A^\Red\times\Gamma(x)': P(\gamma(\infty))<0\right\}&\longrightarrow &\A\langle x\rangle\\
(P,\gamma)&\mapsto &P|\gamma.
\end{array}
$$
\end{teorema}
\textbf{Proof.} We will prove a bijection between the sets on the left hand side in the maps of Theorems \ref{teorema 1} and \ref{teorema 2}. Then Theorem \ref{teorema 2} will follow from Theorem \ref{teorema 1}.

The bijection mentioned above is the map $\psi$:
$$
\begin{array}{rll}
\left\{(P,\gamma)\in\A^\Red\times\Gamma(x)':P(\gamma(\infty))<0\right\}&\stackrel{\psi}{\longrightarrow}&\left\{\begin{array}{ll}(P,\gamma)\in\A^\Sim\times\Gamma(x):\\ P(\gamma(\infty))<0<P(\left\lfloor \gamma(x)\right\rfloor)\end{array}\right\}\\
(P,T^{-k}\gamma_i)&\mapsto&(P|{T^{-k}},\gamma_i).
\end{array}
$$
To prove that $\psi$ is well defined we only have to check, given $(P,\gamma)=\psi(R,\tilde\gamma)$ with $\tilde\gamma=T^{-k}\gamma_i$, two conditions: $P\in\mathcal{A}^\Sim$ and $P(\left\lfloor \gamma(x)\right\rfloor)>0$. We have
$$
w_P=w_R+k,\qquad\qquad w_P'=w_R'+k.
$$
From $w_R'>-1$ and $k\geq 1$, we deduce that $w_P'>0$. The inequality $R(\tilde\gamma(\infty))<0$ and the equality $\tilde\gamma(\infty)=-\dfrac{q_{i-2}}{q_{i-1}}-k$ imply $-w_R>\dfrac{q_{i-2}}{q_{i-1}}+k$, where each term is positive, so $k<-w_R$, and thus $w_P<0$. Hence $P\in\mathcal{A}^\Sim$.
The condition $P(\left\lfloor \gamma(x)\right\rfloor)>0$ follows from
$$
\left\lfloor \gamma(x)\right\rfloor=\left\lfloor \dfrac{\delta_{i-1}}{\delta_i}\right\rfloor=n_{i}\geq k=w'_P-w'_R>w'_P.
$$ 

We consider the map $\varphi$
$$
\begin{array}{rll}
\left\{\begin{array}{l}(P,\gamma)\in\A^\Sim\times\Gamma(x):\\ P(\gamma(\infty))<0<P(\left\lfloor \gamma(x)\right\rfloor)\end{array}\right\}&\stackrel{\varphi}{\rightarrow} &\left\{(P,\gamma)\in\A^\Red\times\Gamma(x)':P(\gamma(\infty))<0\right\}\\
\left(P,\gamma_i\right)&\mapsto&\left(P|{T^{\left\lfloor w'_P\right\rfloor+1}},T^{-\left\lfloor w'_P\right\rfloor-1}\gamma_i\right).
\end{array}
$$
To prove that $\varphi$ is well defined we will check, given $(R,\tilde\gamma)=\varphi(P,\gamma_i)$, two conditions: $R\in\A^\Red$ and $\left\lfloor w'_P\right\rfloor+1\leq n_{i}$. The condition $R\in\A^\Red$ is immediate: from equalities
$$
w_R=w_P-\left\lfloor w_P'\right\rfloor-1,\qquad\qquad w_R'=w_P'-\left\lfloor w_P'\right\rfloor-1
$$
and inequalities $w_P<0<w_P'$, we deduce $w_R<-1<w_R'<0$.

The condition $\left\lfloor w'_P\right\rfloor+1\leq n_{i}$ follows from $\left\lfloor \gamma_i(x)\right\rfloor=n_{i}$ and $P(\left\lfloor \gamma_i(x)\right\rfloor)>0$.

Finally $\varphi$ is the inverse of $\psi$. Indeed, it is clear that $\psi\circ\varphi$ is the identity, and we deduce the same statement for $\varphi\circ\psi$ from the equation below for $(P,\gamma_i)=\psi(R,T^{-k}\gamma_i)$:
$$
-k+\left\lfloor w_P'\right\rfloor+1=-k+\left\lfloor w_R'+k\right\rfloor+1=\left\lfloor w_R'\right\rfloor+1=0.
$$
\begin{flushright}
$\square$
\end{flushright}


\begin{coro}\label{teo igualdades 1}
Let $x$ be a real number, $k\geq 2$ an integer and $\A$ a $\Gamma$-equivalence class in $\mathcal{Q}$. Then the following equalities hold
$$
A_{k,\A}(x)=\sum_{Q\in\mathcal{A}^{\Sim}}\sum_{\substack{\gamma\in\Gamma(x) \\Q(\left\lfloor\gamma(x)\right\rfloor)>0\\Q(\gamma(\infty))<0}}(Q|\gamma)(x)^{k-1}=\sum_{Q\in\A^{\Red}}\sum_{\substack{\gamma\in\Gamma(x)' \\Q(\gamma(\infty))<0}}(Q|\gamma)(x)^{k-1}.
$$
\end{coro}

\begin{coro} For $x\in\mathbb{R}$, the sum $A_{k,D}(x)$ has exponential convergence. It is finite if and only if $x\in\mathbb{Q}$.
\end{coro}


\textbf{Proof.} 
The function $A_{k,D}(x)$ is the sum of the sums that appear in each side of Corollary \ref{teo igualdades 1} over all $\Gamma$-equivalence classes in $\mathcal{Q}_D$. We can prove its exponential convergence looking at the sum on $\A^{\Sim}$ or on $\A^{\Red}$. Let us look at the sum on $\A^{\Sim}$.
On the one hand, the set of polynomials that belong to $\Sim$ with fixed positive discriminant is finite. On the other hand, for an element $Q(X)=aX^2+bX+c$ of $\Sim$, $Q|\gamma_{i}(x)=a\delta_{i-1}^2+b\delta_{i-1}\delta_i+c\delta_i^2$. Now the series $\delta_i=|p_{i-1}-q_{i-1}x|$ stops if $x\in\mathbb{Q}$ and decreases exponentially to 0 if $x\not\in\mathbb{Q}$.
\begin{flushright}
$\square$
\end{flushright}

\begin{teorema}\label{teorema 1b}
For $d\in 2\mathbb{N}$, $\B$ a $\Gamma_1$-equivalence class in $\mathcal{F}_d$ and $x\in\mathbb{R}$, the following bijection holds 
$$
\begin{array}{rll}
\left\{\begin{array}{c}
(P,\gamma)\in\B^\Sim\times\Gamma(x):\\ P(\gamma(\infty))<0<P(\gamma(x))\end{array}
\right\} &\stackrel{\cong}{\longrightarrow} &\B\langle x\rangle\\
(P,\gamma) &\mapsto &P|\gamma.
\end{array}
$$
\end{teorema}

\textbf{Proof.} The map is well defined because of its definition. Let $P\in\B$ satisfy $P(\infty)<0<P(x)$. For $j\gg 0$, the convergents $\dfrac{\tilde{p}_j}{\tilde{q}_j}$ of $x$ belong to $(w_P',w_P)$, because $x\in(w_P',w_P)$. Moreover, if one convergent belongs to $(w_P',w_P)$, so do all the later convergents. 

If $\left\lceil  x\right\rceil\not\in(w_P',w_P)$, we define $i$ to be the unique positive integer such that $\dfrac{\tilde p_{i-1}}{\tilde q_{i-1}}\not\in(w_P',w_P)$ but $\dfrac{\tilde p_i}{\tilde q_i}\in(w_P',w_P)$. If $\left\lceil  x\right\rceil\in(w_P',w_P)$, we set $i=0$. Since $P(\infty)<0$, in both cases we have
\begin{equation}\label{signos}
P\Big(\dfrac{\tilde p_{i-1}}{\tilde q_{i-1}}\Big)<0<P\Big(\dfrac{\tilde p_i}{\tilde q_i}\Big).
\end{equation} 
Put 
$$
\gamma=\begin{pmatrix}-\tilde q_{i-1} &\tilde p_{i-1}\\-\tilde q_i &\tilde p_i\end{pmatrix},
$$
and $R=P|\gamma^{-1}$. We have
$$
R(0)=P\Big(\dfrac{\tilde p_{i-1}}{\tilde q_{i-1}}\Big)<0,\qquad R(\infty)=P\Big(\dfrac{\tilde p_i}{\tilde q_i}\Big)>0,
$$
thus $(R,\gamma)$ belongs on the left hand set of the map in Theorem \ref{teorema 1b}.
The uniqueness of the preimage $(R,\gamma)$ follows from the uniqueness of $i$ satisfying \eqref{signos}.
\begin{flushright}
$\square$
\end{flushright}

\begin{coro}\label{teo igualdades 1b}
Let $x$ be a real number, $k\geq 2$ an integer and $\B$ a $\Gamma_1$-equivalence class in $\mathcal{Q}$. Then
$$
\sum_{Q\in\B\langle x\rangle}Q(x)^{k-1}=\sum_{Q\in\B^\Sim}\sum_{\substack{\gamma\in\Gamma_1(x) \\Q(\gamma(x))>0\\Q(\gamma(\infty))<0}}(Q|\gamma)(x)^{k-1}.
$$
\end{coro}

\section{Continued fractions and modular forms}


There is a canonical choice for the Eichler integral $F$ of a cusp form $f(\tau)=\sum_{n=1}^\infty a_n q^n$ ($q=e^{2\pi i\tau}$) of weight $2k$ ($k\geq 1$) for $\Gamma_1$:
\begin{equation}\label{F canonica}
F(\tau)=\int_\tau^\infty f(z)(\tau-z)^{2k-2}\, dz\, \dot{=}\, \sum_{n=1}^\infty \dfrac{a_n}{n^{2k-1}}\, q^n \qquad(\tau\in\mathcal{H})
\end{equation}
where the symbol $\dot{=}$ denotes equality with constant. The integral defining $F(\tau)$ converges also for $\tau\in\mathbb{R}$, so we can expand the definition domain of $F$ to $\mathcal{H}\cup\mathbb{R}$. The image $F|_{2-2k}(1-\gamma)$ belongs to the space $\mathbb{C}[X]_{\leq 2k-2}$ because of Bol identity between the ($k-1$)st derivative of $F|_{2-2k}\gamma$ and the image by $|_{2k}\gamma$ of the ($k-1$)st derivative of $F$. Moreover, the map
$$
\begin{array}{rll}
\Gamma_1 &\longrightarrow &\mathbb{C}[X]_{\leq 2k-2}\\
\gamma &\mapsto &F|(1-\gamma)
\end{array}
$$ 
is a parabolic 1-cocyle. Since $\Gamma_1$ is generated by $T$ and $S$, the 1-cocyle above is determined by $T\mapsto 0$ and
\begin{equation}\label{rf}
S\, \mapsto\,  r_f(X)=F|(1-S)=\int_0^\infty f(z)(X-z)^{2k-2}\, dz.
\end{equation}
The polynomial $r_f(X)$ is called \textit{period polynomial} of the cusp form $f$. More generally, a map $\Gamma_1\rightarrow\mathbb{C}[X]_{\leq 2k-2}$ which sends $T$ to 0 and $S$ to a complex polynomial $P(X)$ is a parabolic 1-cocyle if and only if $P$ satisfies (\cite{Z00})
\begin{equation}\label{cociclo}
P|(1+S)\=0,\qquad P|(1+U+U^2)\=0.
\end{equation}
We can easily check that if $P=A|(1-S)$, with $A(x)$ a periodic real function, then $P$ satisfies \eqref{cociclo}.

The space $W_{2k}$ of polynomials in $\mathbb{C}[X]_{\leq 2k}$ satisfying \eqref{cociclo} splits up into the subspaces of even and odd polynomials $W^+_{2k}$ and $W^-_{2k}$. Thus  $r_{f}=r_{f}^{+}+r_f^{-}$ with $r_f^{+}(X)\in W_{2k-2}^{+}$ and $r_f^{-}(X)\in W_{2k-2}^{-}$; such polynomials give rise to the isomorphisms
$$
\begin{array}{rrll}
r^{-}: &S_{2k}(\Gamma_1) &\longrightarrow &W_{2k-2}^{-},\\
&f &\mapsto &r_f^{-}(X)
\end{array}
\quad 
\begin{array}{rrll}
r^{+}: &S_{2k}(\Gamma_1) &\longrightarrow &W_{2k-2}^{+}/\left\langle X^{2k-2}-1\right\rangle\\ 
&f &\mapsto &r_f^{+}(X) \pmod{(X^{2k-2}-1)},
\end{array}
$$
where $S_{2k}(\Gamma_1)$ is the space of cusp forms of weight $2k$ for $\Gamma_1$ ($k\geq 1$).
 
In fact we can find the even and odd parts of $r_f$ from the even and odd parts of the Eichler integral of $f$ on $\mathbb{R}$:
\begin{equation}
F^+(x)=\sum^\infty_{n=1}\dfrac{a_n}{n^{2k-1}}\, \cos(2\pi nx),\qquad F^-(x)=\sum^\infty_{n=1}\frac{a_n}{n^{2k-1}}\, \sin(2\pi nx),
\end{equation}
\begin{equation}
r^+_f\equiv F^+|(1-S)\pmod{(X^{2k-2}-1)},\qquad r^-_f=F^-|(1-S).
\end{equation}

Given $P(X)\in\mathbb{C}[X]_{\leq 2k}$, we denote
$$
P^\Gamma(x)\=\sum_{\gamma\in\Gamma(x)}(P|\gamma)(x),\qquad P^{\Gamma_1}(x)\=\sum_{\gamma\in\Gamma_1(x)}(P|\gamma)(x),\qquad(x\in\mathbb{R}).
$$
We will see that when $P(X)\in W^+_{2k}$, the function $P^\Gamma(x)$ is the even part of the Eichler integral on $\mathbb{R}$ of a cusp form whose even part of the period polynomial is (modulo $X^{2k}-1$) the polynomial $-P(X)$. As a consequence, we obtain that we can drop the conditions 
$S(\gamma(\infty))<0<~S(\left\lfloor \gamma(x)\right\rfloor)$ in Corollary \ref{teo igualdades 1}, and $S(\gamma(\infty))<0<S(\gamma(x))$ in Corollary \ref{teo igualdades 1b}. 


For $k\geq 1$ and $P(X)\in\mathbb{C}[X]_{\leq 2k}$, the theorem below gives the differences between $P^\Gamma$ and its image by each generator of $\Gamma$ in terms of $P$, $P|(1+\varepsilon)$ and $P|(1+U-SU^2)$, the last two vanish when $P$ belongs to $W^+_{2k}$.

\begin{teorema}\label{P} For $k\geq 1$ and $P(X)\in\mathbb{C}[X]_{\leq 2k}$, we have
\begin{enumerate}[(i)] 
\item $P^\Gamma|(1-T)=0$.

\item If we denote $P_1=P|(1+\varepsilon)$ and $P_2=P|(1+U-SU^2)\varepsilon$, then $P^\Gamma|(1-\sigma)$ is equal to
$$
\sum_{n\in 2\mathbb{Z}}\Big(\chi_{(\frac{n}{2},\frac{n+1}{2}]}\, (-P_1|T^{n/2+1}\sigma+P_2|T^{-n/2})\, +\, \chi_{(\frac{n+1}{2},\frac{n+2}{2})}\, (P_1|T^{-n/2} - P_2|T^{n/2+1}\sigma)\Big).
$$

\item If $x>0$ and $x\neq 1$, then
$$
(P^\Gamma|(1-\varepsilon))(x)\=\chi_{(0,1)}(P|(1+\varepsilon))(x)\, -\, P(x).
$$
\end{enumerate}
\end{teorema}

\textbf{Proof.} Statement (i) follows from the calculation
\begin{align*}
P^\Gamma(x+1)&\=\sum_{\gamma\in\Gamma(x+1)}(P|\gamma)(x+1)\\
&\=\sum_{\gamma\in\Gamma(x)}(P|\gamma T^{-1})(x+1)\qquad\mbox{because $\Gamma(x+1)=\Gamma(x)T^{-1}$}\\
&\=P^\Gamma(x).
\end{align*}

If $x\in\mathbb{Z}$, then $P^\Gamma(x)=P(\infty)$, so $P^\Gamma|(1-\sigma)(x)=0$.

Let $x=[n_0,n_1,\ldots]$ be a non integer number and $-x=[-n_0-1,n_1',\ldots]$ its opposite.  The inequality  $0<x-n_0\leq\frac{1}{2}$ is equivalent to $n_1\geq 2$. It is also equivalent to $\frac{1}{2}\leq-x+n_0+1<1$, so to $n_1'=1$.
By Proposition~\ref{prop 1}, we have
$$
P^\Gamma(-x)=\sum_{\substack{-1\leq\gamma(\infty)\leq0\\\gamma(-x)>1}}(P|\gamma)(-x)\,\, -\,\, (P|\varepsilon T^{-n_0})(x)\,  
-\, \left\{\begin{array}{cc}
(P|T^{-1}\varepsilon T^{n_0+1}\sigma)(x) &\, \mbox{if $n_1=1$}\\
\emptyset &\, \mbox{if $n_1\geq 2$}.
\end{array}\right.
$$
Since
$$
\sum_{\substack{-1\leq\gamma(\infty)\leq0\\\gamma(-x)>1}}(P|\gamma)(-x)\=\sum_{\substack{-1\leq\gamma \sigma(\infty)\leq0\\\gamma \sigma(x)>1}}(P|\gamma \sigma)(x)\=\sum_{\substack{-1\leq\gamma(\infty)\leq0\\\gamma(x)>1}}(P|\gamma)(x),
$$
by Proposition \ref{prop 1} again, we have
$$\begin{array}{lll}
P^\Gamma|(1-\sigma)&\=&P|\varepsilon T^{-n_0}\, -\, P|\varepsilon T^{n_0+1}\sigma
\, +\, \left\{\begin{array}{cc}
-P|T^{-1}\varepsilon T^{-n_0} &\quad\mbox{if $n_1\geq 2$}\\
P|T^{-1}\varepsilon T^{n_0+1}\sigma &\quad\mbox{if $n_1=1$}
\end{array}\right.\\\\
&\=&\left\{\begin{array}{ll}-P|(1+\varepsilon)T^{n_0+1}\sigma\, +\, P|(1+U-SU^2)\varepsilon T^{-n_0} &\quad\mbox{si $n_1\geq 2$}\\
P|(1+\varepsilon)T^{-n_0}\, -\, P|(1+U-SU^2)\varepsilon T^{n_0+1}\sigma &\quad\mbox{si $n_1=1$}
\end{array}\right.\end{array}
$$
because $T^{-n_0}=U\varepsilon T^{n_0+1}\sigma$, $\, T^{n_0+1}\sigma=U\varepsilon T^{-n_0}$ and $SU^2=T^{-1}$. Thus statement (ii) is proved.

Suppose $x>1$. For $i\geq 1$, the $i$-th term of the real series $(x_i)_{i\geq 0}$ defined in \eqref{algo Red 1} which gives the continued fraction of $1/x$ is equal to the $(i-1)$-th term of the series which gives the continued fraction of $x$. So
$$
\Gamma(1/x)\=\Gamma(x)\, \varepsilon\, \cup\, \left\{\begin{pmatrix}0 &1\\1 &0\end{pmatrix}\right\}.
$$
From that we get
\begin{align*}
x^{2k}P^\Gamma(1/x)&\=x^{2k}\sum_{\gamma\in\Gamma(x)}(P|\gamma \varepsilon)(1/x)\, +\, x^{2k}(P|\varepsilon)(1/x)\\
&\=P^\Gamma(x)+P(x).
\end{align*}

Now suppose $0<x<1$ and denote $y=1/x$. We have
\begin{align*}
x^{2k}P^\Gamma(1/x)&\=\dfrac{1}{y^{2k}}\, P^\Gamma(y)\\
&\=P^\Gamma(1/y)-\dfrac{1}{y^{2k}}\, P(y)\qquad\mbox{(by the previous case)}\\
&\=P^\Gamma(x)-x^{2k}P(1/x).
\end{align*}
\begin{flushright}
$\square$
\end{flushright}

\begin{coro}\label{prop funcion auxiliar} 
\begin{enumerate}[(i)]
\item  For $k\geq 1$ and $P(X)\in W^+_{2k}$, the function $P^\Gamma(x)$ is even, periodic and satisfies $P^\Gamma|(1-S)=-P$. 

\item Let $f$ be a cusp form of weight $2k$ for $\Gamma_1$ and $P(X)$ be the even part of its period polynomial. The even part of the Eichler integral of $f$ on $\mathbb{R}$ is (modulo a constant multiple of $X^{2k-2}-1$) the function $(-P)^\Gamma(x)$.
\end{enumerate}
\end{coro}

\textbf{Proof.} (i) Let $P(X)$ be an element of $W^+_{2k}$. By the statement (i) of Theorem \ref{P}, the function $P^\Gamma(x)$ is periodic.
Since $1+U-SU^2=1+U+U^2-(1+S)U^2$ and $P(X)$ is even (so $P|S=P|\varepsilon$), statement (ii) implies that $P^\Gamma(x)$ is even. 

By statement (iii), we have $P^\Gamma|(1-S)(x)=-P(x)$ for $x\neq 0,1$. If $x=0$, then $(P^\Gamma|(1-S))(0)=P(\infty)=-P(0)$ because $P|(1+S)=0$. If $x=1$, then $(P^\Gamma|(1-S))(1)=0=P(1)$ again because $P|(1+S)=0$ and $P|(1-\sigma)=0$.

(ii) Let $f$ be a cusp form of weight $2k$ for $\Gamma_1$ and $P(X)$ be the even part of its period polynomial. Then $P(X)$ belongs to $W^+_{2k-2}$, so $(-P)^\Gamma$ is periodic and satisfies $(-P)^\Gamma|(1-S)=P$ by statement (i) of the corollary. Hence $(-P)^\Gamma(x)$ is (modulo $X^{2k-2}-1$) the even part of the Eichler integral of $f$ for $x\in\mathbb{R}$.
\begin{flushright}
$\square$
\end{flushright}

Given an even integer $k\geq 2$ and a $\Gamma$-equivalence class $\A$  of~ $\mathcal{Q}$, we define the polynomial 
$$
P_{k,\A}(X)=\sum_{Q\in\A^\Sim}Q(X)^{k-1}.
$$
Kohnen and Zagier proved in \cite{KZ84} that for every positive non-square discriminant $D$, the polynomial 
$$
P_{k,D}(X)=\sum_{\A\in\mathcal{Q}_D/\Gamma}P_{k,\A}(X)
$$ 
is (modulo a constant multiple of $X^{2k-2}-1$) the even part of the period polynomial of the cusp form of weight $2k$ for the modular group $\Gamma_1$
$$
f_{k,D}(z)\=\sum_{\substack{(a,b,c)\in\mathbb{Z}^3\\b^2-4ac=D}}\dfrac{1}{(az^2+bz+c)^k}\qquad(z\in\mathcal{H},\, k\geq 2\, \mbox{even}).
$$
This function arose in \cite{Z75}, in the case where $D$ was a fundamental discriminant, considering the restriction to the diagonal $z_1=z_2$ of a family of Hilbert modular forms $w_m(z_1,z_2)$ ($m=0,1,2,\ldots$) of weight $k$ for the Hilbert modular group $\mathrm{SL}_2(\O)$, where $\O$ was the ring of integers of the real quadratic field with discriminant $D$. The functions $w_m(z_1,z_2)$ are the Fourier coefficients of the kernel function for Doi-Naganuma correspondence between elliptic modular forms and Hilbert modular forms. They are well defined for all positive discriminants $D$ and so is $f_{k,D}(z)$, except that when $D$ is a square there is an extra term besides $P_{k,D}$ in the expression of the even part of the period polynomial.

Kohnen and Zagier proved in \cite{KZ80} that the functions $D^{k-1/2}\, f_{k,D}(z)$ are the $D$-th Fourier coefficients of the kernel function for Shimura-Shintani correspondence. 

Recently, Bringmann, Kane and Kohnen gave a new proof of the fact that $P_{k,D}(X)$ is the even part of the period polynomial of $f_{k,D}(z)$ for positive non-square discriminants $D$ using new modular objects related to the theory of harmonic weak Maass forms (\cite{BKK}). 

In \cite{Z99} (§6), Zagier proved that $P_{k,D}(X)$ belongs to $W^+_{2k-2}$ showing that the function $A_{k,D}(x)$ is periodic and satisfies $A_{k,D}|(1-S)=-P_{k,D}$. The same argument he used applies to $A_{k,\A}(x)$ if $\A$ is a $\Gamma$-equivalence class in $\mathcal{Q}$, so $A_{k,\A}|(1-S)=-P_{k,\A}$ and for even $k\geq 2$ the polynomial $-P_{k,\A}(X)$ belongs to $W^+_{2k-2}$. Hence $-P_{k,\A}$ is (modulo $X^{2k-2}-1$) the even part of a period polynomial and $A_{k,\A}(x)$ is the corresponding even part of the Eichler integral on $\mathbb{R}$. Then, by Corollary \ref{prop funcion auxiliar}, $P^\Gamma_{k,\A}(x)=A_{k,\A}(x)$ for all $x\in\mathbb{R}$. Together with Theorem \ref{teo igualdades 1}, we obtain

\begin{coro} For a $\Gamma$-equivalence class $\A$ of $\mathcal{Q}$ and an even integer $k\geq 2$, the following identities hold:
$$
A_{k,\A}(x)\=\sum_{Q\in\A^\Sim}\, \sum_{\gamma\in\Gamma(x)} (Q|\gamma)(x)^{k-1}\=\sum_{Q\in\A^\Sim}\, \sum_{\substack{\gamma\in\Gamma(x)\\Q(\gamma(\infty))<0\\Q(\left\lfloor \gamma(x)\right\rfloor)>0}} (Q|\gamma)(x)^{k-1}.
$$
\end{coro}

\begin{coro} For every positive non-square discriminant $D$, the functions $P^\Gamma_{2,D}(x)$ and $P^\Gamma_{4,D}(x)$ have the respective values $-5L(-1,\chi_D)$ and $L(-3,\chi_D)$.
\end{coro}

\begin{Remarque}\label{Q y S} We can also give a direct proof of the identity $P^\Gamma_{2,D}(x)=-5L(-1,\chi_D)$. To do this, write $P^\Gamma_{2,D}(x)$ as
$$
P^\Gamma_{2,D}(x)\=\sum_{i\geq 1}\sum_{\substack{(a,b,c)\in\mathbb{Z}^3\\a>0>c\\b^2-4ac=D}}(a\delta_{i-1}^2+b\delta_{i-1}\delta_i+c\delta_i^2).
$$
If $(a,b,c)$ appears in the sum, then $(-c,-b,-a)$ appears too. Hence
$$
P^\Gamma_{2,D}(x)\=\Big(\sum^\infty_{i=1}\, (\delta_{i-1}^2-\delta_i^2)\Big)\cdot\Big(\sum_{\substack{a,b,c\in\mathbb{Z}\\a>0>c\\b^2-4ac=D}} a\Big).
$$
But the first sum telescopes to 1 because the $\delta_i$ decrease to 0 and $\delta_{0}=1$, and the second sum equals $-5L(-1,\chi_D)$ by
results of \cite{S69}, \cite{C76} and \cite{Z76}.
\end{Remarque}

If we sum over the quadratic functions now in a single $\Gamma_1$-equivalence class $\B$, we have to symmetrize with respect to an involution on the set $\mathcal{Q}_D/\Gamma_1$ to construct functions related to modular forms. For each $\Gamma_1$-equivalence class $\B$, and each (not necessarily even) integer $k\geq 2$, we define
$$
A_{k,\B}^\ast(x)=A_{k,\B}(x)\, +\, (-1)^k A_{k,-\B}(x),
$$
$$
P_{k,\B}(X)=\sum_{Q\in\B^\Sim}Q(X)^{k-1}\, +\, (-1)^k\sum_{Q\in(-\B)^\Sim}Q(X)^{k-1},
$$
where $-\B:=\left\{-Q\mid\, Q\in\B\right\}$.

Again we can use the argument in \S6 of \cite{Z99} to deduce that $A_{k,\B}^\ast(x)$ is periodic and satisfies $A_{k,\B}^\ast|(1-S)=-P_{k,\B}$ for $k\geq 2$. Thus 
$$
-P_{k,\B}|\tilde\gamma_{i+1}\=(A_{k,\B}^\ast-A_{k,\B}^\ast|T^{m_i}S)|\tilde\gamma_{i+1}\=A_{k,\B}^\ast|\tilde\gamma_{i+1}-A_{k,\B}^\ast|\tilde\gamma_i\qquad (i\geq 0)
$$
and so $A_{k,\B}^\ast(x)=P^{\Gamma_1}_{k,\B}$. Hence we obtain

\begin{coro} For a $\Gamma_1$-equivalence class $\B$ of $\mathcal{Q}$ and an integer $k\geq 2$, the following identities hold
\begin{align*}
A_{k,\B}^\ast(x)&\=\sum_{Q\in\B^\Sim}\, \sum_{\gamma\in\Gamma_1(x)} (Q|\gamma)(x)^{k-1}\, +\, (-1)^k\sum_{Q\in(-\B)^\Sim}\, \sum_{\gamma\in\Gamma_1(x)} (Q|\gamma)(x)^{k-1}\\
&\=\sum_{Q\in\B^\Sim}\, \sum_{\substack{\gamma\in\Gamma_1(x)\\Q(\gamma(\infty))<0\\Q(\gamma(x))>0}} (Q|\gamma)(x)^{k-1}\, +\, (-1)^k\sum_{Q\in(-\B)^\Sim}\, \sum_{\substack{\gamma\in\Gamma_1(x)\\Q(\gamma(\infty))<0\\Q(\gamma(x))>0}} (Q|\gamma)(x)^{k-1}.
\end{align*}
\end{coro}

\section{Remarks on the quartic case}

We could try to give a similar construction to the function $A_{k,D}(x)$ for sums taken over quartic polynomials. Quartic polynomials can have 0, 2 or 4 real roots; in the second case the discriminant is negative, otherwise it is positive. The analog case to the quadratic construction is to consider polynomials with at least two real roots, because for such a polynomial $Q$ satisfying $Q(\infty)<0$, the set of real numbers $x$ on which $Q(x)\geq 0$ is compact, as for the quadratic case. 

The discriminant can be written in terms of the two $\Gamma_1$-invariants $I$ and $J$ associated to a quartic polynomial $aX^4+bX^3+cX^2+dX+e$:
$$
I=12ae-3bd+c^2,\qquad J=72ace+9bcd-27ad^2-27eb^2-2c^3,\qquad D=\dfrac{1}{27}(4I^3-J^2).
$$

A naive generalization of the function $A_{k,D}(x)$ would be taking the sum over the $k$-th powers of quartic polynomials with integer coefficients, fixed $I$ and $J$, that are negative at $\infty$ and positive at $x$. But this sum diverge. In fact by Theorem \ref{teorema 1}, for a $\Gamma$-equivalence class $\A$ in $\mathcal{F}_4$ and an even integer $k\geq 1$,
\begin{equation}\label{cuartico}
\sum_{P\in\A\left\langle x\right\rangle} P(x)^k\=\sum_{P\in\A^{\Sim}}\, \sum_{\substack{\gamma\in\Gamma(x)\\P(\left\lfloor \gamma(x)\right\rfloor)>0\\P(\gamma(\infty))<0}}(P|\gamma)(x)^k\=\sum_{P\in\A^{\Red}}\,  \sum_{\substack{\gamma\in\Gamma(x)'\\P(\gamma(\infty))<0}}(P|\gamma)(x)^k
\end{equation}
where neither the set $\A^{\Sim}$ nor $\A^{\Red}$ is finite. 

We could try to modify the naive generalization modifying the left or the right hand side of equation \eqref{cuartico}. On the right hand side, we should replace $\A^{\Sim}$ by a finite set $\A^{\mathrm{Fin}}$ of polynomials in $\A$ such that $\sum_{P\in\A^{\mathrm{Fin}}}P(X)^k$ is the even part of a period polynomial. Because $\sum_{P\in\A^{\mathrm{Fin}}}P(X)^k$ should be invariant by $1+S$, the power $k$ should be odd (as for the quadratic case).

If we look at the left hand side, we should add some linear inequalities for the coefficients of the polynomials $P(X)$ to make the sum converge. But any linear inequality involving other coefficients of $P(X)$ than $P(\infty)$ would probably break the invariance by $T$ of the sum, because the only invariants by $T$ for a quartic $P(X)$ are $P(\infty)$, $I$, $J$ and $P(X)$ itself. So the new sum would not be the even part of an Eichler integral anymore.


\end{document}